\documentclass[11pt]{amsart}
\usepackage{epsfig,graphics,amsmath,amssymb,amsthm}
\usepackage{graphics}

\newtheorem{theorem}{Theorem}[section]
\newtheorem{lemma}[theorem]{Lemma}
\newtheorem{proposition}[theorem]{Proposition}
\newtheorem{corollary}[theorem]{Corollary}

\theoremstyle{definition}

\newtheorem{remark}[theorem]{Remark}

\theoremstyle{remark}

 \newcommand{\abs}[1]{\lvert#1\rvert}
 \def\R{{\mathbb{R}}}
 \def\Z{{\mathbb{Z}}}

 \def\S{{\mathbb{S}}}
 \def\mod{{\rm Mod}}

\begin{document}

\newenvironment{prooff}{\medskip \par \noindent {\it Proof}\ }{\hfill
$\square$ \medskip \par}
    \def\sqr#1#2{{\vcenter{\hrule height.#2pt
        \hbox{\vrule width.#2pt height#1pt \kern#1pt
            \vrule width.#2pt}\hrule height.#2pt}}}
    \def\square{\mathchoice\sqr67\sqr67\sqr{2.1}6\sqr{1.5}6}
\def\pf#1{\medskip \par \noindent {\it #1.}\ }
\def\endpf{\hfill $\square$ \medskip \par}
\def\demo#1{\medskip \par \noindent {\it #1.}\ }
\def\enddemo{\medskip \par}
\def\qed{~\hfill$\square$}


 \title[Lefschetz fibrations and finitely presented groups]
 {Lefschetz fibrations and an invariant of finitely presented groups}

 \author{Mustafa Korkmaz}

 \address{Department of Mathematics, Middle East Technical University,
 06531 Ankara, Turkey} \email{korkmaz@metu.edu.tr}

\thanks{The author is supported in part by the Turkish Academy of Sciences
under the Young Scientists Award Program (MK/T\"UBA-GEB\.IP 2003-10).}
 \date{\today}

\begin{abstract}
Every finitely presented group is the fundamental group of the total space of a 
Lefschetz fibration. This follows from results of Gompf and Donaldson, and was also 
proved by Amoros-Bogomolov-Katzarkov-Pantev. We give another proof by providing 
the monodromy explicitly. We then define the genus of a finitely presented group 
$\Gamma$ to be the minimal genus of a Lefschetz fibration with fundamental group
$\Gamma$. We also give some estimates of the genus of certain groups.
  
\end{abstract}
 \maketitle
 \setcounter{secnumdepth}{2}
 \setcounter{section}{0}

\section{Introduction} \label{section1}
Every finitely presented group is the fundamental group of some closed symplectic 
$4$-manifold~\cite{g}. By the work of Donaldson~\cite{d}, every closed oriented 
symplectic $4$-manifold admits a Lefschetz pencil. Thus such a manifold admits 
a Lefschetz fibration over $\S^2$, perhaps after blowing up many times. It
follows that any finitely presented group is realized as the fundamental group
of the total space of a Lefschetz fibration over $\S^2$. Conversely, if
$g\geq 2$, then the total space of a genus $g$ Lefschetz fibration is symplectic (c.f.~\cite{gs}).

The paper~\cite{abkp} gives another 
construction of a symplectic $4$-manifold as the total space 
of a Lefschetz fibration with the given fundamental group. In this construction,
the genus of the Lefschetz fibration depends on the number of intersection points
of some curves on some surface representing the relators of the group. Thus it is quadratic in the 
lengths of the relators.

The purpose of this paper is to give yet another construction of a Lefschetz fibration 
with the prescribed finitely presented group. In our construction, we give the monodromy 
explicitly.  The genus of the Lefschetz fibration depends linearly on the number of generators
and on the syllable lengths of the relators of the presentation of the given finitely presented group.
More precisely, given a finitely presented group $\Gamma$ with $n$ generators and $k$ relators,
if $\ell$ is the sum of the syllable lengths of the relators, then for every 
$g\geq 2(n+\ell -k)$ we construct a genus-$g$ Lefschetz fibration
over $\S^2$ whose fundamental group is isomorphic to $\Gamma$. 

\begin{figure}[hbt]
\label{figure:1}
 \begin{center}
    \includegraphics[width=6cm]{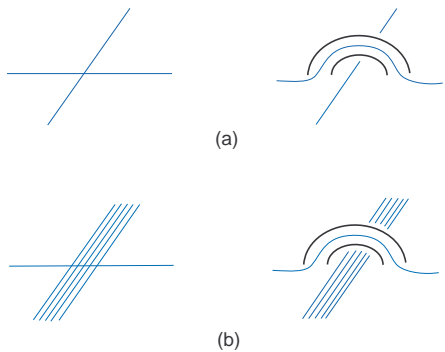}
  \caption{Two ways of resolving the self intersection points of a loop.}
  \label{operation}
   \end{center}
 \end{figure}

In~\cite{abkp}, the authors first represent the relators of the finitely presented group
by loops on a certain surface $\Sigma$ and then for each intersection point of 
these loops they increase the genus of $\Sigma$ by one as illustrated in Figure~\ref{figure:1}~(a). 
This makes the genus of the fiber of the Lefschetz fibration 
quadratic in the lengths of the relators. However, in our proof,
we increase the genus of $\Sigma$ by one for each set of consecutive self intersection points
of a loop as illustrated in Figure~\ref{figure:1}~(b).

Here is an outline of the paper. In Section~\ref{section2}, we fix the notations, state the main theorem,
Theorem~\ref{thm:main}, 
and give the necessary tool from the theory of mapping class groups
to prove it. Section~\ref{section3} reviews the theory of 
Lefschetz fibrations and proves Theorem~\ref{thm:presentation}, which is
used throughout the paper. Section~\ref{section4} is devoted to the proof of
Theorem~\ref{thm:main}. In Section~\ref{section5}, 
we define an invariant of a finitely presented group $\Gamma$,
which we call the genus of $\Gamma$. It turns out that the genus of the fundamental group
of a closed orientable surface of genus $g$ is $g$. We compute the genera of the groups
$\Z, \Z\oplus \Z_n, \Z_m \oplus \Z_n, \Z_n$, which are all equal to $2$. We also give upper 
and lower bounds for the genus for any finitely presented group and we obtain some results 
of the genera of some groups related to topology; free groups, finitely generated abelian groups, 
fundamental groups of closed nonorienrable surfcaes, braid groups 
and the group $SL(z,Z)$. In the final section, Section~\ref{section6}, we compare our 
invariant with Kotschick's invariants $p$ and $q$, and pose some problems motivated by the
content of the paper. 

I would like to thank to Andras Stipsicz and Burak Ozbagci for the conversations regarding the 
content of this paper.

\section{Preliminaries and the statement of main result}
\label{section2}
Let $\Gamma$ be a finitely generated group generated by a set
$X=\{ x_1,x_2,\ldots,x_n\}$. For an element $w\in \Gamma$, let the
\textit{syllable length}  $\ell(w)$ of $w$ be defined as follows:
$$
\ell(w)=\min \{ \, s\, |\,\, w=x_{i_1}^{m_1}x_{i_2}^{m_2}\cdots
 x_{i_s}^{m_s}, \, 1\leq i_j\leq n, \, m_j\in\Z  \}. $$

Suppose now that $\Gamma$ is a finitely presented group with a
presentation
\begin{eqnarray}
\Gamma=\langle x_1,x_2,\ldots,x_n | \,\, r_1,r_2,\ldots, r_k
\rangle, \label{eqn:pres}
\end{eqnarray}
so that $\Gamma$ is the quotient $F/N$, where $F$ is the free
group (nonabelian for $n\geq 2$) freely generated by $\{ x_1,\ldots, x_n\}$ and $N$ is the
normal subgroup of $F$ generated normally by the elements
$r_1,r_2,\ldots, r_k$. That is, $N$ is the subgroup of $F$ generated by
all conjugates of $r_1,r_2,\ldots, r_k$.

 Define $\ell=\ell(r_1)+\ell(r_2)+\cdots +\ell(r_k)$, which depends on
the presentation. We always assume that the relators $r_i$ are
cyclically reduced.

The main result of this paper is the following theorem.

\begin{theorem} \label{thm:main}
Let $\Gamma$ be a finitely presented group with a presentation
$(\ref{eqn:pres})$. Then for every $g\geq 2(n+\ell-k)$ there exists a
genus-$g$ Lefschetz fibration $f:X\to \S^2$ such that
$\pi_1(X)$ is isomorphic to $\Gamma$. 
\end{theorem}

For a closed oriented surface $\Sigma_g$ of genus $g$, the mapping class group
$\mod_g$ of $\Sigma_g$ is defined to be the group of isotopy classes of orientation
preserving diffeomorphism $\Sigma_g\rightarrow \Sigma_g$. If $a$ is a simple closed 
curve on $\Sigma_g$, then cutting $\Sigma_g$ along $a$ and gluing the two boundary components
back after twisting one of the sides to the right by $2\pi$ give a diffeomorphism $\Sigma_g \rightarrow \Sigma_g$.
The isotopy class of this diffeomorphism, denoted by $t_a$, is called the right 
Dehn twist about $a$. The mapping class $t_a$ depends only on the isotopy class of $a$.
Dehn twists are the simplest diffeomorphisms of $\Sigma_g$. They are the main building blocks
in $\mod_g$; they generate the mapping class group and each representation of the identity element
$1$ in the group $\mod_g$ as a product of right Dehn twists gives a Lefschetz fibration, whose
total space is a symplectic $4$-manifold. 

We now describe such a relation in $\mod_g$ which is the main ingredient in the 
proof of our result. This relation was obtained in~\cite{k} as an extension of Matsumoto's
relation in~\cite{m} and was used to construct noncomplex smooth symplectic $4$-manifolds 
admitting genus-$g$ Lefschetz fibrations. It was also used by Ozbagci and Stipsicz in~\cite{os}
in their construction of infinitely many Stein fillings of a certain contact $3$-manifold.

\begin{figure}[hbt]
 \begin{center}
    \includegraphics[width=8cm]{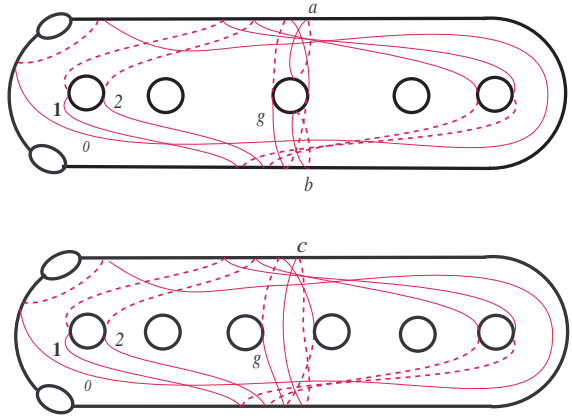}
  \caption{The simple closed curve labelled $j$ is $B_j$.}
  \label{figure1}
   \end{center}
 \end{figure}

\begin{figure}[hbt]
 \begin{center}
    \includegraphics[width=9cm]{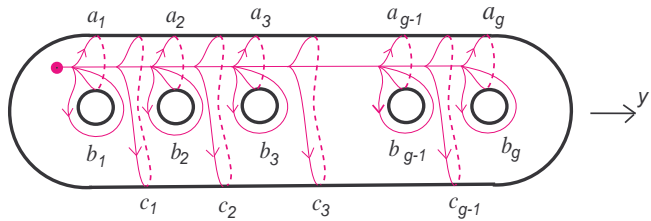}
  \caption{Generators of the fundamental group.}
  \label{figure2}
   \end{center}
 \end{figure}

Let $\Sigma_g$ denote the closed oriented surface of genus $g$
standardly embedded in the $3$-space as in Figures~\ref{figure2}, 
so that it is the boundary a $3$-dimensional handlebody.
Let $\Sigma_{g,2}$ denote the surface of genus $g$ with two boundary components
obtained from $\Sigma_g$ by deleting the interior of two disjoint discs (c.f. 
Figure~\ref{figure1}). Let us define a word
 \begin{eqnarray} \label{wrodW}
 W=\left\{
 \begin{array}{ll}
 (t_c^2t_{B_g}t_{B_{g-1}}\cdots t_{B_2} t_{B_1}t_{B_0})^2 & {\rm if}\,\,g\,\,
 {\rm is\,\,even,}\\
 (t_a^2t_b^2 t_{B_g}t_{B_{g-1}}\cdots t_{B_2} t_{B_1}t_{B_0})^2 & {\rm if}\,\,g\,\,
 {\rm is\,\, odd},
 \end{array}
 \right.
 \end{eqnarray}
in the mapping class group of $\Sigma_{g,2}$, 
where the simple closed curves $B_j$ and $a,b,c$ are given in Figure~\ref{figure1}. 
It was shown in~\cite{k} that the word $W$ represents the identity element
in the mapping class group of $\Sigma_g$. It can be shown easily that 
the word $W$ is equal to the product $t_{\delta_1}t_{\delta_2}$ in the 
mapping class group of $\Sigma_{g,2}$. Here, $\delta_1$ and $\delta_2$
are the boundary components of $\Sigma_{g,2}$. When there is one boundary
component $\delta$, it was shown in~\cite{os2} that the word $W$ is equal to 
the Dehn twist $t_{\delta}$, and this will be sufficient for us.

We note that Dehn twists in the relation $W$ above is given in the reversed order
in~\cite{k}. There the functional notation was used for the composition of functions. 
However, in this paper, the composition 
$fg$ of two diffeomorphism $f$ and $g$ means that we first apply $f$ and then $g$. 

Consider the standard set of generators $a_1,b_1,a_2,b_2,\ldots, a_g,b_g$ of
$\pi_1(\Sigma_g)$ as shown in Figure~\ref{figure2}.
The subgroup of $\pi_1(\Sigma_g)$ generated by $a_1,a_2,\ldots, a_g$
is a free group of rank $g$.

Throughout the paper a loop and its homotopy class will be denoted by the same notation.
Similarly, a diffeomorphism and its isotopy class or a simple closed curve and its isotopy
class will be denoted by the same symbol. We even denote a simple loop and a simple closed curve by
the same symbol, and this will not cause any problem as it will be clear from the context
which one we mean.

If $x$ and $y$ are two elements of a group, we write $[x,y]=x^{-1}y^{-1}xy$ and 
$x^y=y^{-1}xy$. We recall that if $f$ is in mapping class group of an oriented surface 
$\Sigma$ and $t_a$ is a Dehn twist about a simple closed curve $a$ on $\Sigma$, 
then it easily follows from the definition of a Dehn twist
that $f^{-1}t_af=t_{f(a)}$, so that for a mapping class (in particular for a Dehn twist) $f$,
the conjugation $W^{f}=1$ is a representation of the identity as a product of right Dehn twists
in Mod$_g$.

\section{Lefschetz fibrations}
\label{section3}
Let us review briefly the theory of Lefschetz fibrations on $4$-manifolds.
For the details the reader is referred to~\cite{gs}.
Let $X$ be a closed connected oriented smooth $4$-manifold. A
Lefschetz fibration on $X$ is a smooth map $f:X\to \S^2$ such that
if $q$ is a critical point, then $q$ has a neighborhood
$U$ with complex coordinates $(z_1,z_2)$ agreeing with the orientation of 
$X$ and $f(q)$ has a neighborhood $V$ with complex coordinate $z$ agreeing with the orientation of 
$\S^2$ such that the restriction $f|:U\to V$ of $f$ is of the form
$f(z_1,z_2)=z_1^2+z_2^2$.

It follows that the number of critical points of $f$ is finite. We
can assume that each singular fiber contains only one critical
point. Suppose that $Q=\{ q_1,q_2,\ldots, q_s\}$ are critical
points and $f(Q)=\{ p_1,p_2,\ldots, p_s\}$ are the critical values
so that $f(q_i)=p_i$. Removing the singular fibers from the total
space $X$ and the critical values from the base, we get a surface
bundle 
$$X-{f^{-1}(f(Q))}\longrightarrow \S^2- {f(Q)}.$$ 
Therefore, all regular
fibers of $f$ are diffeomorphic to a closed connected oriented
smooth surface of genus $g$ for some $g$. We will assume that all
Lefschetz fibrations are relatively minimal, i.e. no fiber
contains an embedded $2$-sphere with self-intersection number
$-1$. 

Let us fix a base point $p_0$ on $\S^2$ which is a regular
value and let $\Sigma$ denote the fiber over it; $\Sigma=f^{-1}(p_0)$. If $C$ is a
simple loop on $\S^2-f(Q)$ encircling only one of the critical
values, then $f^{-1}(C)$ is surface bundle over a circle with the
monodromy a right Dehn twist about a simple closed curve on
$\Sigma$, which is called a vanishing cycle.

The topology of the total space of a Lefschetz fibration is determined 
by the monodromy representation $\varphi:\pi_1(\S^2-f(Q))\to \mod_g$.
For a simple loop $C$ encircling only one critical value,
if $c$ is the corresponding vanishing cycle,
then $\varphi (C)$ is the Dehn twist $t_c$. 

Let $D$ be a disc on $\S^2$ containing all critical values in the
interior and the base point $p_0$ on the boundary. For each $i=1,2,\ldots, s$  
let $e_i$ be a small circle in the interior of $D$ 
encircling the critical value $p_i$ but no other critical value. 
Let $d_i$ be a simple arc connecting $p_0$ to a point $e_i$.
We assume that $d_i$ and $d_j$ intersect only at $p_0$ for $i\neq j$.
We choose $d_i$ in such a way that, on the boundary $\partial$ of a small 
regular neighborhood of $p_0$, the intersection of $\partial$ and 
$d_1,d_2,\ldots,d_s$ are read in the counter clockwise order.
Now let $\alpha_i$ be the loop obtained first travelling from $p_0$
to $e_i$ along $d_i$, then travelling along $e_i$ once counter clockwise
and then back to $p_0$ along $d_i$. Clearly, the fundamental group
$\pi_1(\S^2-f(Q))$ has a presentation with generators 
$\alpha_1,\alpha_2,\ldots, \alpha_s$ and with a single defining relation 
$$\alpha_1 \alpha_2\cdots \alpha_s=1.$$ Thus 
if $c_i$ is the vanishing cycle corresponding to $\alpha_i$, then the Dehn twists
$t_{c_i}$ satisfy the relation 
 \begin{eqnarray}\label{prod=1}
 t_{c_1}t_{c_2}\cdots t_{c_s}=1
 \end{eqnarray}
in the mapping class group $\mod_g$ of the surface $\Sigma$.

Conversely, given a decomposition of the identity as a product of right Dehn  
twists of the form~(\ref{prod=1}) in the mapping class group $\mod_g$,
we can construct a genus-$g$
Lefschetz fibration as follows: Consider the product
$\Sigma_g\times D$, where $\Sigma_g$ is a closed connected
oriented surface of genus $g$ and $D$ is a $2$-disc. For each
$i=1,2,\ldots,s$, attach a $2$-handle to $\Sigma_g\times D$ along
the simple closed curve $c_i$ with framing $-1$ relative to the
product framing in $\Sigma_g\times \partial D$. The resulting
$4$-manifold $Y$ admits a Lefschetz fibration over the $2$-disc
whose monodromy around the boundary is the product
$t_{c_1}t_{c_2}\cdots t_{c_s}$, which is equal to the identity in
the mapping class group of $\Sigma_g$. Therefore, the boundary of
$Y$ is diffeomorphic to the product $\Sigma_g\times S^1$. Now glue $\Sigma_g\times
D$ and $Y$ along their boundaries to get the closed connected
oriented smooth $4$-manifold $X$ admitting a Lefschetz fibration
with generic fiber $\Sigma_g$.

A section of a Lefschetz fibration $f:X\to \S^2$ is a map $\sigma :\S^2\to X$
such that $f\sigma$ is the identity map of $\S^2$.  Since the word $W$ represents the 
element $t_{\delta_1}t_{\delta_2}$ in the mapping class group of $\Sigma_{g,2}$,
the Lefschetz fibration of genus $g$ with the monodromy $W=1$ admits two disjoint 
sections of self intersection $-1$ (c.f.~\cite{smith}).

\begin{lemma}
\label{lemma:pi1} $($ \cite{gs} $)$
Let $f:X\to \S^2$ be a genus-$g$ Lefschetz fibration with global
monodromy given by the relation~$(\ref{prod=1})$. Suppose that
$f$ has a section. Then the
fundamental group of $X$ is isomorphic to the fundamental group of
$\Sigma_g$ divided out by the normal closure of the simple closed
curves $c_1,c_2,\ldots,c_s$, considered as elements in
$\pi_1(\Sigma_g)$. In particular, there is an epimorphism 
$\pi_1(\Sigma_g)\to \pi_1(X)$
 
\end{lemma}

Suppose that $V$ is a product of right Dehn twists representing the identity
in the mapping class group Mod$_g$ of the closed oriented surface $\Sigma_g$. 
If $f_1,f_2,\ldots, f_m$ are arbitrary elements in Mod$_g$, then the product
$$V^{f_1}V^{f_2}\cdots V^{f_m}$$ 
is also a product of right Dehn twists representing the identity in Mod$_g$. Here,
$V^{f}$ denotes the conjugation $f^{-1} V f$. Let us denote 
by $X^g_V(f_1,f_2,\ldots,f_m)$ the total space of the Lefschetz fibration over $\S^2$ 
whose monodromy is $V^{f_1}V^{f_2}\cdots V^{f_m}$. For simplicity,
if $d_i$ is a simple closed curve on $\Sigma$, we denote the $4$-manifold
$X^g_V(t_{d_1},t_{d_2},\ldots,t_{d_m})$ by $X^g_V(d_1,d_2,\ldots,d_m)$.
We note that $X^g_V(d_1,d_2,\ldots,d_m)$ is a fiber sum 
$$X^g_V(1)\#_fX^g_V(1)\#_f \cdots \#_fX^g_V(1)$$
of $m$ copies of $X^g_V(1)$. 
Clearly, if $X^g_V(1)$ has a section, then $X^g_V(d_1,d_2,\ldots,d_m)$
has a section as well. 

The main result of this section is the following theorem, which is used throuhout the paper.

\begin{theorem}
\label{thm:presentation}
Let $V= t_{c_1}t_{c_2}\cdots t_{c_n}$ be a product of right Dehn 
twists representing the identity in the mapping class group 
of the closed oriented surface $\Sigma_g$. Suppose that
the Lefschetz fibration $X^g_V(1)\to \S^2$ has a section. 
Let $d_1,d_2,\ldots,d_m$ be simple closed curves on $\Sigma_g$ with the 
property that for each $j=1,2,\ldots,m$, there exists an $i_j$ such that 
$d_j$ intersects $c_{i_j}$ transversely
at only one point. Then the fundamental group of $X^g_V(1,d_1,d_2,\ldots,d_m)$
is isomorphic to the group $\pi_1(X^g_V(1))$ divided out by 
the normal closure of $\{ d_1,d_2,\cdots,d_m\}$.
\end{theorem}

\begin{proof}
By Lemma~\ref{lemma:pi1}, the group $\pi_1(X^g_V(1))$ admits a presentation with generators
$a_1,b_1,\ldots,a_g,b_g$ and with relations
 \begin{itemize}
  \item $[a_1,b_1]^{y_1}\cdots [a_g,b_g]^{y_g}=1$ and
  \item $c_1=c_2=\cdots =c_n=1$,
\end{itemize}
where $y_i$ are some words in $a_j,b_j$. Since $V^f= t_{f(c_1)}t_{f(c_2)}\cdots t_{f(c_n)}$, 
in order to obtain a presentation 
of the group $\pi_1(X^g_V(1,d_1,d_2,\ldots,d_m))$ one needs to add the relations $t_{d_j} (c_i)=1$
to this presentation. That is, 
the fundamental group $\pi_1(X^g_V(1,d_1,d_2,\ldots,d_m))$ admits a presentation with generators
$a_1,b_1,\ldots,a_g,b_g$ and with relations
\begin{eqnarray} \label{eqn:4}
  && [a_1,b_1]^{y_1}\cdots [a_g,b_g]^{y_g}=1, \nonumber\\ 
  && c_1=c_2=\cdots =c_n=1, \mbox{ and} \\
  && t_{d_j} (c_i)=1, \,\, 1\leq i\leq n, \, 1\leq j\leq m. \nonumber
\end{eqnarray}

Let us fix some $1\leq j\leq m$. By assumption, the simple closed curve $d_j$ intersects 
$c_l$ transversely at one point for some $1\leq l\leq n$. Fix an orientation of $d_j$, 
Then the curve $t_{d_j} (c_l)$ is equal (in fact conjugate) to $c_l\, d_j^\varepsilon$, where $\varepsilon ={\pm 1}$. 
Since $c_l=1$ in this presentation, we may replace the relator $t_{d_j} (c_l)=1$ by 
$d_j=1$. Suppose that $i\neq l$. If $d_j$ is disjoint from $c_i$, then $t_{d_j} (c_l)=c_l$.
Suppose now that $d_j$ intersects $c_i$ at $t$ points. It can be seen easily that 
there are elements $x_1,x_2,\ldots,x_{t+1}$ in $\pi_1 (\Sigma_g)$ such that 
$c_i=x_1x_2\ldots x_{t+1}$ and that the curve $t_{d_j} (c_l)$ is conjugate 
to $x_1d_j^{\varepsilon_1}x_2d_j^{\varepsilon_2}\ldots x_td_j^{\varepsilon_t} x_{t+1}$,
where each $\varepsilon_i$ is equal to $1$ or $-1$. Since $d_j=1$ and $c_i=1$, we can delete the
relators $t_{d_j}(c_i)=1$ from the presentation (\ref{eqn:4}). 

Repeating this argument for each $j=1,2,\ldots,m$ shows that the relators  
$t_{d_j} (c_i)=1$, $1\leq i\leq n$, $1\leq j\leq m$, may be replaced by 
$d_j=1$. 

This proves the theorem.
\end{proof}

\section{The proof of Theorem~\ref{thm:main} }
\label{section4}
In this section, we prove Theorem~\ref{thm:main}.
For the proof, we take appropriate fiber sums of the 
Lefschetz fibration $X^g_W(1)$ with global monodromy $W$ with itself.

\subsection{The proof of Theorem~\ref{thm:main} for free groups}
From the presentation point of view, the simplest groups are free groups.
So we assume first that $\Gamma$ is a free group.

\begin{proposition} \label{prop}
Let $\Gamma$ be a free group of rank $n$. Then for
every $g\geq 2n$ there is a genus-$g$ Lefschetz fibration $f:X\to
\S^2$ with $\pi_1(X)\cong\Gamma$.
\end{proposition}

\begin{proof} 
Assume first that $g=2r$ is even, so that $r\geq n$. Let $a_1,b_1,\ldots,a_g,b_g$ be
the standard generators of the fundamental group $\pi_1(\Sigma_g)$ 
of the surface $\Sigma_g$ as shown in Figure~\ref{figure2}.  
It can easily be shown that in $\pi_1(\Sigma_g)$, up to conjugation, the
following equalities hold:
 \begin{itemize}
  \item $B_0=b_1b_2\cdots b_g,$
  \item $B_{2k-1}=a_{k}b_kb_{k+1}\cdots b_{g+1-k}
  c_{g+1-k}a_{g+1-k}$, $1\leq k\leq r,$
  \item $B_{2k}=a_{k}b_{k+1}b_{k+2}\cdots b_{g-k}
  c_{g-k}a_{g+1-k}$, $1\leq k\leq r,$
  \item $c=c_r$.
  \end{itemize}
Note that $c_k=[a_1,b_1]^{y_1}\cdots [a_k,b_k]^{y_k}$ for some
elements $y_1,\ldots,y_k$ in the group $\pi_1(\Sigma_g)$. 

Let us consider the symplectic $4$-manifold 
$$X^g_W(1,b_{1},b_{2},\ldots,b_{g}, a_{n+1},a_{n+2},\ldots,a_{r} )$$
and let us denote it simply by $X$. Notice that for each $i=1,2,\ldots ,g$, the simple closed curve $b_i$ 
intersects two of the curves $B_{j}$ transversely at only one point
(and is disjoint from the others). Similarly, the curve $a_i$ intersects at least one $B_j$
transversely at one point.

By Theorem~\ref{thm:presentation}, the fundamental group of $X$ 
admits a presentation with generators $a_1,b_1,\ldots, a_g,b_g$ and with relations
 \begin{itemize}
  \item $[a_1,b_1]^{y_1}[a_2,b_2]^{y_2}\cdots [a_g,b_g]^{y_g}=1,$
  \item $B_0=B_1=\cdots = B_{g}=c=1,$
  \item $b_{1}=b_{2}=\cdots =b_{g}=1,$ and
  \item $a_{n+1}=a_{n+2}=\cdots =a_{r}=1,$.
 \end{itemize}
It is easy to show that this is a presentation 
of the free group of rank $n$, with a free basis $a_1,a_2,\ldots,a_n$. Thus 
the group $\pi_1(X)$ is a free group of rank $n$.

Suppose now that $g=2r+1$ is odd and $g>2n$. A similar computation shows that 
the fundamental group of 
$$X^g_W(1,b_{1},b_{2},\ldots,b_{g}, a_{n+1},a_{n+2},\ldots,a_{r} )$$
is, again, a free group of rank $n$, finishing the proof of the proposition.
\end{proof}

This proves Theorem~\ref{thm:main} if $\Gamma$ is a free group.

\bigskip

\remark
The number of singular fibers of the Lefschetz fibration in above proposition 
is 
\[
(2g+4)(1+g+(\frac{g}{2}-n))=(2g+4)(1+3r-n)\] 
if $g=2r$, and  
\[
(2g+10)(1+g+(\frac{g-1}{2}-n))=(2g+10)(2+3r-n)\] if $g=2r+1$.
It is possible to reduce the number of singular fibers. In fact, one 
can show that the fundamental group of the $4$-manifold
$$X^g_W(1,b_{r+1},b_{r+2},\ldots,b_{g}, a_{n+1},a_{n+2},\ldots,a_{r} )$$
is also isomorphic to $F_n$.

\subsection{The proof of Theorem~\ref{thm:main} for arbitrary finitely presented group}
Recall that the elements $a_1,a_2,\ldots ,a_n$ of $\pi_1(\Sigma_n)$
shown in Figure~\ref{figure2} generate a free subgroup of rank $n$. In fact, any subset of 
the standard generators with cardinality less than $2n$ generate
a free subgroup. 

Let us begin with the following proposition, which is a part of the proof 
but it deserves to be stated separately, as it is interesting in itself.

\begin{proposition}
\label{prop:R}
Let $F_n$ denote the subgroup of $\pi_1(\Sigma_n)$ generated by 
the generators $a_1,a_2,\ldots ,a_n$, so that $F_n$ is a free group of rank $n$.
Let $r_1, r_2,\ldots, r_k$ be arbitrary $k$ elements in $F_n$
represented as words in $a_1,a_2,\ldots ,a_n$. Let
$\ell$ denote $\ell(r_1)+\ell(r_2)+\cdots +\ell(r_k)$, the sum of the 
syllable lengths of $r_i$, and let $h=n+\ell -k$. Then,
on the closed connected orientable surface $\Sigma_h$ of genus $h$, 
there are loops $R_1, R_2,\ldots, R_k$ 
satisfying the following conditions:
\begin{itemize}
\item[(a)] Each $R_i$ is a simple loop on $\Sigma_h$.

\item[(b)] Each $R_i$ is freely homotopic to a simple closed
curve intersecting $a_h$ transversely at exactly one point.

\item[(c)] If $\Phi:\pi_1(\Sigma_h)\to \pi_1(\Sigma_n)$ denotes the map defined by
$\Phi(a_j)=a_j$ for $1\leq j\leq n$ and $\Phi(\alpha)=1$ for 
$\alpha\in \{a_{n+1},\ldots,a_h, b_1,\ldots,b_h \}$, then 
$\Phi ([R_i])=r_i$ for each $i$, where $[R_i]\in \pi_1(\Sigma_h)$ is the 
homotopy class of $R_i$.
\end{itemize}
\end{proposition}

\begin{figure}[hbt]
 \begin{center}
    \includegraphics[width=9cm]{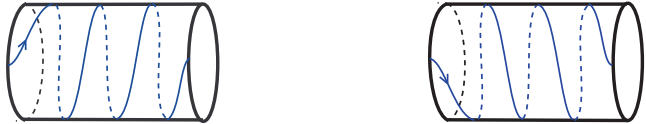}
  \caption{Annuli $P(3)$ and $P(-3)$.}
  \label{annulus}
   \end{center}
 \end{figure}

\begin{proof}
For an integer $m$, let $P(m)$ denote the annulus
$[0,1]\times S^1$ together with the arc $\gamma_m(t)=(t,e^{-2m\pi i t})$,
$0\leq t\leq 1$, on it.

Let us consider the surface $\Sigma_n$ embedded in $\R^3$ as in Figure~\ref{figure2}.
Note that for each $i=1,2,\ldots,n$ there is a constant $\nu_i$ such that
the intersection of $\Sigma_n$ with the plane $y=\nu_i$ is
the disjoint union of two oriented simple closed curves, one of which is
freely homotopic to the loop $a_i$.  Let us denote these simple 
closed curves by $\alpha_i$ and $\alpha'_i$ so that $\alpha_i$
is freely homotopic to $a_i$ and $\alpha'_i$ is obtained from $\alpha_i$
by rotating the surface $\Sigma_n$ by $\pi$ about the $y$-axis.
Note that in general we denote $\alpha_i$ by $a_i$; we will use distinct notations
only in this proof.

Let $L$ be an embedded arc on $\Sigma_n$ connecting the base point
to a point on $\alpha_n$ such that $L$ intersects each $\alpha_i$ at
a single point and is disjoint from each $\alpha'_i$. We
can assume that the arc $L$ lies in a plane $z$=constant.

Now let $r=a_{i_1}^{m_1} a_{i_2}^{m_2}\cdots a_{i_d}^{m_d}$ be an element 
of the free group $F_n$, where $d=\ell(r)$ is the syllable length of $r$. 
Choose pairwise disjoint simple closed curves
$\tilde{\alpha}_{i_1},\tilde{\alpha}_{i_2},\ldots,\tilde{\alpha}_{i_d}$
such that $\tilde{\alpha}_{i_j}$ is isotopic to $\alpha_{i_j}$. 
We can assume further that
$\tilde{\alpha}_{i_j}$ is the intersection of a plane 
$y=\nu_{i_j}$ with the surface $\Sigma_n$, and if $i_j=i_k$ for some
$j<k$, then $\nu_{i_j}<\nu_{i_k}$. Moreover, we assume that each
$\tilde{\alpha}_{i_j}$ intersects $L$ at exactly one point.

In what follows, we glue various copies of the annulus $P(m)$ to
$\Sigma_n$ minus the interiors of disjoint open discs and annuli
in such a way that the resulting surface is orientable, and the obvious
orientation of $P(m)$ and that of $\Sigma_n$ agree.

For each $j=1,2,\ldots, d$, let us choose a regular neighborhood
$N_j$ of $\tilde{\alpha}_{i_j}$ so that the closures of any two
distinct $N_j$ are pairwise disjoint and do not contain the base
point. We now identify $N_j$ and the annulus $P(m_{j})$ so that
$[0,1]\times  1$ is identified with a subarc of $L$ and the $y$-coordinate of
$0\times 1$ is less than that of $1\times 1$. This gives us a bunch
of disjoint arcs on $\Sigma_n$ as shown in Figure~\ref{construction}(a).

On the annulus $N_j$, let us label the starting point
of the curve $\gamma_{m_{j}}$ by $A_j$ and the terminal
point by $B_j$. We note that the $y-$coordinate of $A_j$
is less than that of $B_j$. For each $1\leq j\leq d-1$,
let $\delta_j$ denote the subarc of $L$ from the point $B_j$ 
to the point $A_{j+1}$. Then
$$
\beta=\gamma_{m_{1}}\star\delta_1\star\gamma_{m_{2}}\star\delta_2\star
\cdots \star\delta_{d-1}\star\gamma_{m_{d}}
$$
is an arc on $\Sigma_n$ connecting $A_1$ to $B_d$, and if we delete
all $\delta_j$ from $\beta$, the result is the disjoint union of $d$ simple arcs.
Furthermore, if $\delta_0$ is the subarc of $L$ from the base point to 
$A_1$ and $\delta_d$ is the subarc of $L$ from $B_d$ to the base point,
then the loop $\delta_0\star \beta\star \delta_d$ is a representative of $r$.  
By perturbing $\beta$ slightly, we assume that 
$\delta_0,\delta_1,\delta_2,\ldots,\delta_d$ are pairwise disjoint.

It is now clear from our construction that the base point can be joined
to the point $A_1$ by an arc $\delta'$ and the point $B_d$ can be joined
to the base point by an arc $\delta''$
such that interiors of $\beta,\delta'$ and $\delta''$ are pairwise disjoint and
that the loop $\delta'\star \beta\star\delta''$ represents the element
$$
b_1b_2\cdots b_{i_1-1}\, r \,  b_{i_d}^{-1}\cdots b_2^{-1} b_1^{-1}
$$
in $\pi_1(\Sigma_n)$ (c.f. Figure~\ref{construction}(b)). We can assume that the 
arc $\delta'$ intersects $\alpha'_1, \ldots,\alpha'_{i_1-1}$ at one point
and is disjoint from all other $\alpha'$ and all $\alpha$.
Similarly, the arc $\delta''$ intersects only $\alpha'_1, \ldots,\alpha'_{i_d}$
at one point.

Let $D_1,D_2,\ldots,D_{2d-2}$ be pairwise disjoint discs on $\Sigma_n$
such that the interior Int$(D_i)$ of each $D_i$ is disjoint from $\beta$,
$\delta'$ and $\delta''$,
and for each $j=1,2,\ldots,d-1$, $\partial D_{2j-1}\cap \beta=\{ B_j\}$
and $\partial D_{2j}\cap\beta =\{ A_{j+1}\}$.

\begin{figure}[hbt]
 \begin{center}
    \includegraphics[width=13cm]{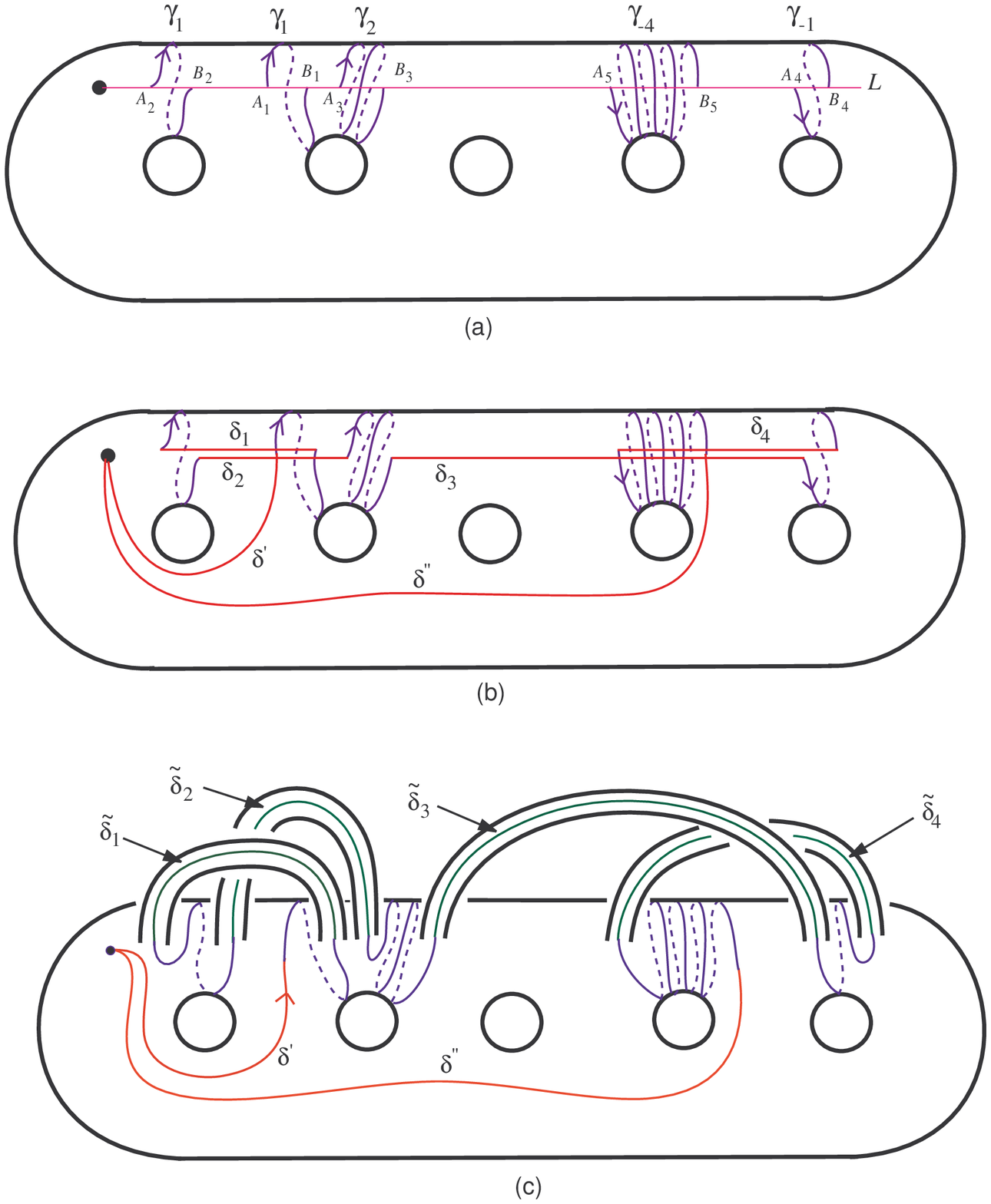}
  \caption{Construction of $R$ on $\Sigma_{n+d-1}$ for $r=x_2x_1x_2^2x_5^{-1}x_4^{-4}$ and 
for $n=5$.}
  \label{construction}
   \end{center}
 \end{figure}

We now remove $2d-2$ open discs Int$(D_i)$ from $\Sigma_n$ and
for each $1\leq j\leq d-1$ we glue a copy of the annulus $P(0)$ to the surface
$$
\Sigma_n\backslash \bigcup_{i=1}^{2d-2} \mbox{ Int}(D_i)
$$
in such a way that $0\times S^1$ (respectively $1\times S^1$) is identified 
with $\partial D_{2j-1}$ (respectively $\partial D_{2j}$) and
the point $0\times 1$ (respectively $1\times 1$) is identified with $B_j$ 
(respectively $A_{j+1}$). In this way, we obtain the closed orientable surface
$$
\left( \Sigma_n\backslash \bigcup_{i=1}^{2d-2} \mbox{ Int}(D_i)\right)\cup
 \left( \bigcup_{j=1}^{d-1} P(0) \right)
$$
of genus $n+d-1$, which is denoted by $\Sigma_{n+d-1}$. The orientation on 
$\Sigma_n$ gives an orientation on $\Sigma_{n+d-1}$.

Let $R$ be the loop on $\Sigma_{n+d-1}$ obtained as follows.
For each $j=1,2,\ldots, d-1$, let $\tilde{\delta}_j$  denote the arc $\gamma_0=I\times 1$
on the $j$th annulus $P(0)$, so that it is a simple arc joining the point
$B_j$ to $A_{j+1}$ on $\Sigma_{n+d-1}$. Then the loop
$$
R=\delta'\star \gamma_{m_{1}}\star
 \tilde{\delta}_1\star\gamma_{m_{2}}\star\tilde{\delta}_2\star
\cdots \star\tilde{\delta}_{d-1}\star\gamma_{m_{d}}\star\delta''
$$
obtained from $\delta'\star\beta\star\delta''$ by ``replacing" $\delta_j$ by $\tilde{\delta}_j$
is simple on $\Sigma_{n+d-1}$ (c.f. Figure~\ref{construction}(c)). 

We note that it follows from the construction that $\tilde{\delta}_j\star\delta_j^{-1}$
is a simple closed curve on $\Sigma_{n+d-1}$.
Collapsing each $P(0)$ onto the arc $\delta_j$ gives
a map $\Sigma_{n+d-1}\to \Sigma_n$ and the induced map
between the fundamental groups takes $[R]$ to 
$b_1b_2\cdots b_{i_1-1}\, r \,  b_{i_d}^{-1}\cdots b_2^{-1} b_1^{-1}$,
which is mapped to $r$ under $\Pi$.

Clearly, the above construction can be done for all elements
$r_i$ simultaneously instead of a single element of $F_n$. 
For each $i=1,2,\ldots,k$, we increase the genus of the surface by $\ell (r_i)-1$.
Thus the resulting surface is a closed orientable surface $\Sigma_{h}$
of $h=n+\ell -k$. We have $k$ simple loops $R_1,R_2,\ldots, R_k$
on $\Sigma_h$ such that the homomorphism $\Phi:\pi_1(\Sigma_{h})\to \pi_1(\Sigma_n)$ 
maps $[R_i]$ to $r_i$ for all $i$.

Now slide the $\ell -k$ cylinders that we attached on $\Sigma_n$  
to bring the surface into the standard position as shown in Figure~\ref{figure2}
 so that for each $j=1,2,\ldots,\ell-k$
the simple closed curve $\delta'_j{\delta_j}^{-1}$ becomes isotopic to 
$b_{n+j}$ and the core $\frac{1}{2}\times S^1$ of the $j$th handle 
becomes isotopic to $a_{n+j}$. Note that exactly one of $R_i$ intersects 
the curve $a_h$ transversely only once. It is clear that those $R_i$ that does not intersect
$a_h$ can be modified to intersect $a_h$ at one point at the expense of 
multiplying $[R_i]$ by some elements $b_j$ for $j>n$, which are mapped 
to the identity under $\Phi$.   

This finishes the proof of the proposition.
\end{proof}

\subsection*{Finishing the proof}

We now continue with the proof of Theorem~\ref{thm:main}. Let $g\geq 2h$ be an integer,
where $h=n+\ell-k$.

Suppose first that $g$ is even, say $g=2r$. Suppose that for each $i=1,2,\ldots,k$, 
the relator $r_i$ in the presentation (\ref{eqn:pres}) is represented by 
the word $V_i(x_1,x_2,\ldots,x_n)$. Let us denote the word $V_i(a_1,a_2,\ldots,a_n)$
also by $r_i$.

Consider the surface $\Sigma_h$ and the loops
$R_j$ constructed in Proposition~\ref{prop:R} using these $r_i$. 
Let us remove the interior of a small disc from $\Sigma_h$ 
near the curve $a_h$ and disjoint from all $R_j$. 
Denote by $\Sigma_{h,1}$ the resulting surface of genus $h$ with one boundary component.
Embed $\Sigma_{h,1}$ into our standard surface $\Sigma_g$ so that simple loops
$a_1,a_2,\ldots ,a_h,b_1,b_2,\ldots ,b_h$ on the surface $\Sigma_{h,1}$ with one boundary 
component coincide with the loop on $\Sigma_g$ labelled by the same letters.
Thus, the loops $R_j$ on $\Sigma_g$ are disjoint from the simple closed curves
$a_{h+1},b_{h+1},\ldots, a_{g},b_{g}$.

It follows that, on the surface $\Sigma_g$, the simple closed curve $B_{2h}$ intersects each $R_j$ 
transversely once.

We note that the element $[R_i]\in \pi_1(\Sigma_{g})$ is contained in the subgroup
generated by $a_{1},\ldots,a_r$ and $b_{1},\ldots,b_r$. Moreover, if we set $a_s=1$ 
for $n+1\leq s\leq r$ and $b_j=1$ for $1\leq j\leq r$
in a word representing $[R_i]\in \pi_1(\Sigma_{g})$, we get a word representing 
the element $r_i$.

\begin{figure}[hbt]
 \begin{center}
    \includegraphics[width=14cm]{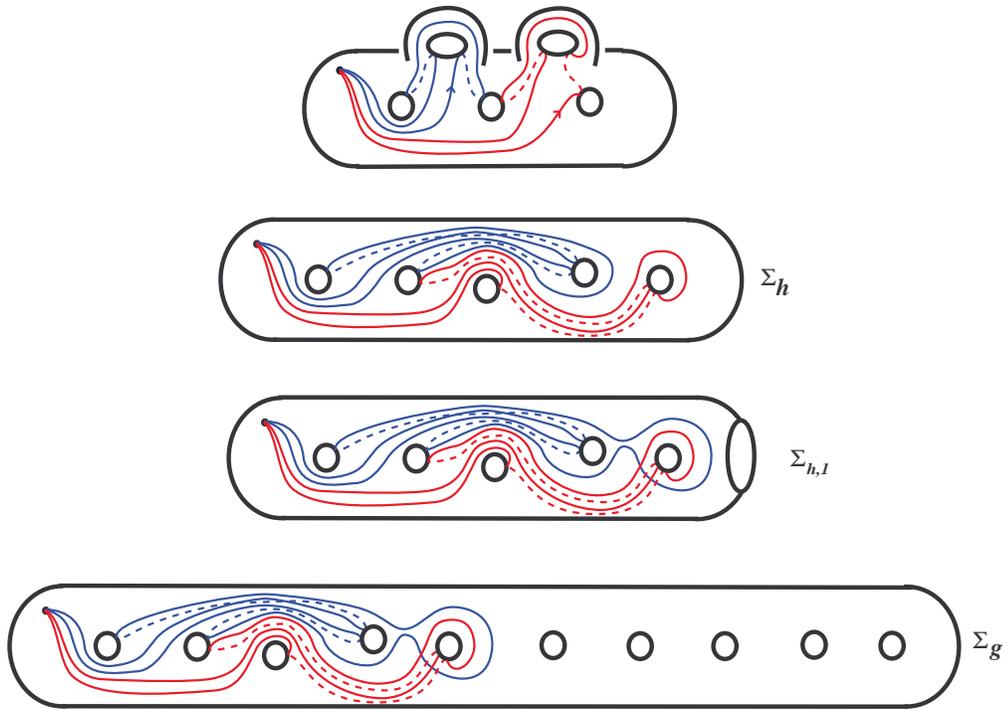}
  \caption{Construction of $R_1$ and $R_2$ for $r_1=x_2x_1^{-1}$ and $r_2=x_3^{-1}x_2^{-1}$
in the case $n=3$ and $g=10$.}
  \label{kaydirma}
   \end{center}
 \end{figure}

Let $X$ denote the $4$-manifold 
$$X^g_W(1,b_{1},\ldots,b_{g}, a_{n+1},\ldots,a_{r},R_1\ldots,R_k )$$
which admits the structure of a Lefschetz fibration. Since for each of
the simple closed curves $b_{1},\ldots,b_{g}, a_{n+1},\ldots,a_{r},R_1\ldots,R_k$
there is at least one $B_j$ such that they intersect precisely at one point, 
by Theorem~\ref{thm:presentation} $\pi_1(X)$ admits a presentation with generators
$a_1,\ldots,a_g, b_1,\ldots,b_g$ and with the defining relations
 \begin{itemize}
  \item $[a_1,b_1]^{y_1}[a_2,b_2]^{y_2}\cdots [a_g,b_g]^{y_g}=1;$
  \item $B_0=B_1=\cdots = B_{g}=c=1;$
  \item $b_i=1$ for $ 1\leq i \leq g;$
  \item $a_i=1$ for $ n+1\leq i \leq r;$
  \item $R_i=1$ for $ 1\leq i \leq k$.
 \end{itemize}
This presentation is equivalent to the presentation~(\ref{eqn:pres})
of $\Gamma$.

It can also be shown similarly that the fundamental group of 
$$X^g_W(1,b_{1},\ldots,b_{g}, a_{n+1},\ldots,a_{r},R_1\ldots,R_k )$$
is again isomorphic to $\Gamma$ if $g$ is odd and greater than
$h$.

This completes the proof of Theorem~\ref{thm:main}.

\remark It is clear from the above proof that, our bound $2(n+\ell-k)$ is not optimal;
depending on the presentation, the genus of the Lefschetz fibration can be made smaller. 
We give an example of this sort when we given an upper bound for the genus of
finitely generated abelian groups in Theorem~\ref{thm:freeabelian} and of the 
braid group in the next section.

\remark 
For each positive integer $m$, let $X_m$ denote the $4$-manifold 
$$X^g_W(1,1,\ldots ,1,b_{1},\ldots,b_{g}, a_{n+1},\ldots,a_{r},R_1\ldots,R_k ),$$
where $1$ is repeated $m$ times. Then the fundamental group 
of $X_m$ is isomorphic to $\Gamma$ for all $m$. 
Moreover, if $m\neq l$ then $X_m$ is not diffeomorphic to
$X_l$, because, for example, their Euler characteristics are different.

\bigskip

\section{An invariant of finitely presented groups} \label{section5}
Using the fact that every finitely presented group is the fundamental group of the 
total space of a Lefschetz fibration, we define an invariant of finitely presented groups.
We then compute the invariant of some groups and give bounds for certain groups; free groups,
abelian groups, braid groups and and the group $SL(2\Z)$. 

\subsection{Definition and the genus of certain groups}
In this section, we consider only those Lefschetz fibrations 
which have sections. For a finitely presented group $\Gamma$, 
we define the {\it genus} $g(\Gamma)$ of
$\Gamma$ to be the minimum $g$ such that there exists
a genus-$g$ Lefschetz fibration $X$ over $\S^2$ such that $\pi_1(X)$ is
isomorphic to $\Gamma$.  Since for each finitely presented group $\Gamma$ 
there is a Lefschetz fibration $f:X\to\S^2$ with a section such that $\pi_1(X)\cong \Gamma$, 
the number $g(\Gamma)$ always exists and is a nonnegative integer.
We note that in the definition, we allow a Lefschetz fibration not
to have any singular fibers. In that case, the total space is just the
product of the fiber surface with the base $\S^2$. 

\begin{theorem} \label{thm:surfacegp}
Let $\pi_g$ denote the fundamental group of a closed orientable surface 
of genus $g$. Then
\begin{enumerate}
	\item[(a)] $g(\pi_g)=g$.
	\item[(b)] $g(\Gamma)=0$ if and only if $\Gamma$ is the trivial group.
	\item[(c)] $g(\Gamma)=1$ if and only if $\Gamma$ is isomorphic to $\Z\oplus\Z$.
	\item[(d)] $g(\Gamma)=2$ if $\Gamma\in \{ \Z, \Z\oplus \Z_n, 
	\Z_m\oplus \Z_n, \Z_n\}$, where $m$ and $n$ are positive integers and $\Z_k$ denotes
	the group $\Z / k \Z$.
\end{enumerate}
\end{theorem}

\begin{proof}
 If $X\to\S^2$ is a genus-$h$ Lefschetz fibration $X\to \S^2$ with
$\pi_1(X)=\pi_g$, then there is an epimorphism $\pi_h \to
\pi_g$. Since no map $\pi_h \to \pi_g$ can be surjective unless $h\geq g$, 
it follows that $g\leq g(\pi_g)$. On the other hand, since the
projection $\Sigma_g\times \S^2\to \S^2$ is a genus-$g$ Lefschetz
fibration, we conclude that $g(\pi_g)=g$. This proves (a).

In particular, the genus of the trivial group is zero. Suppose that
$g(\Gamma)=0$. Let $X\to \S^2$ be a genus-$0$ Lefschetz fibration
with $\pi_1(X)=\Gamma$. Then $\pi_1(X)$ is trivial, proving (b).

Let $X\to \S ^2$ be a genus-$1$ Lefschetz fibration. 
If there is no singular fibers then $X=T^2\times \S^2$. If there 
are singular fibers then $X=E(n)$ for some $n$. Hence $X$ is simply connected.
It follows that $g(\Gamma)=1$ if and only if $\Gamma=\Z\oplus \Z$.

Thus, if $\Gamma$ is a nontrivial group other than $\Z\oplus \Z$, then 
$g(\Gamma)\geq 2$. On the other hand, there are genus-$2$ Lefschetz fibrations with 
fundamental groups $\Z\oplus \Z_n$, $\Z$, $\Z_m\oplus \Z_n$ and $\Z_n$. More precisely, 
the Lefschetz fibrations $X^2_W(1,t_{b_1}^n)$, $X^2_W (1,t_{b_1})$,
$X^2_W(1,t_{a_1}^m,t_{b_1}^n)$ and $X^2_W(1,t_{a_1},t_{b_1}^n)$ have fundamental 
groups isomorphic to $\Z\oplus\Z_n$,  $\Z$, $\Z_m\oplus\Z_n$ and $\Z_n$,
respectively. The fact that the fundamental group of $X^2_W(1,t_{b_1}^n)$ is isomorphic to
$\Z\oplus\Z_n$ was shown in~\cite{os}. This proves (d).
\end{proof}

\bigskip
\subsection{Upper and lower bounds}
Let $\Gamma$ be a finitely presented group with a given presentation with $n$ generators
and with $k$ relators. Suppose that $\ell$ is the sum of the syllable lengths
of the relators. Set $d(\Gamma)=n+\ell - k$. The integer $d(\Gamma)$
depends on the presentation. 

\begin{theorem} \label{thm:genus}
Let $\Gamma$ be a nontrivial finitely presented group. 
Let $m(\Gamma)$ denote the minimal number of generators
for $\Gamma$. Then
\begin{enumerate}
	\item[(a)] $g(\Gamma)\leq \inf \{ 2d(\Gamma) \}$, where the infimum is 
	taken over all presentations of $\Gamma$.
	\item[(b)] $\frac{m(\Gamma)}{2} \leq g(\Gamma)$, with the equality
	if and only if $\Gamma$ is a surface group.
\end{enumerate}
\end{theorem}

\begin{proof}
By Theorem~\ref{thm:main}, for a given presentation pf $\Gamma$ there is a Lefschetz fibration 
whose fundamental group is isomorphic to $\Gamma$ and whose fiber 
genus is $2d(\Gamma)$. Thus, $g(\Gamma)\leq 2d(\Gamma)$ for any finite presentation 
of $\Gamma$, proving (a).

Let $X\to\S^2$ be a Lefschetz fibration of genus $h=g(\Gamma)$ with 
$\pi_1(X)\cong \Gamma$. Then there is an epimorphism $\pi_1(\Sigma_h)\to \Gamma$.
Since $\pi_1(\Sigma_h)$ is generated by $2h$ elements, it follows that
minimal number of generator of $\Gamma$ satisfies $m(\Gamma)\leq 2h$.

If $\Gamma=\pi_1 (\Sigma_h)$ is a surface group, then the minimal number of 
generators of $\Gamma$ is $m(\Gamma)=2h=2g(\Gamma)$.
If $\Gamma$ is not a surface group, then the Lefschetz fibration $X\to\S^2$
of genus $h=g(\Gamma)$ has singular fibers. Since the monodromy group of any 
Lefschetz fibration cannot be contained in the Torelli group
 by the work of Smith~\cite{smith}, there must be at least one nonseparating 
vanishing cycle. It follows that $\Gamma$ can be generated by $2h-1$ elements,
that is, $m(\Gamma)\leq 2g(\Gamma)-1$.
\end{proof}

\begin{corollary} 
Let $\Gamma$ be a nontrivial finitely presented group. 
Let $b_1(\Gamma)$ denote the first Betti number of $\Gamma$, the dimension of the 
first homology of $\Gamma$ with rational coefficients. Then
	$\frac{b_1(\Gamma)}{2} \leq g(\Gamma)$.
\end{corollary}

\bigskip
\subsection{Free groups} As usual, let $F_n$ denote the free group
of rank $n$ freely generated by $x_1,x_2,\ldots,x_n$. 

Next theorem is due to Zieschang~\cite{z}.

\begin{theorem} $($\cite{z}$)$\label{thm:ontoFn}
If there is an epimorphism $\varphi:\pi_1(\Sigma_g)\to F_n$, then $n\leq g$.
\end{theorem}
\begin{proof}
The induced map $\varphi^* :H^1(F_n;\Z) \to H^1(\pi_1(\Sigma_g);\Z)$ is injective.
Since the cup product is zero in $H^1(F_n;\Z)$, its image is Lagrangian in 
$H^1(\pi_1(\Sigma_g);\Z)$. Thus $n\leq g$.
\end{proof}

\begin{theorem}
Let $F_n$ denote the nonabelian free group of rank $n$. Then 
$n\leq g(F_n)\leq 2n$. 
\end{theorem}

\begin{proof}
 The group $F_n$ has a presentation with $n$ generators and with no relations. 
Hence, $d(F_n)=n$ for this presentation. Thus, $g(F_n)\leq 2n$  by Theorem~\ref{thm:genus}~(a). 
On the other hand, if there is an epimorphism from $\pi_1 (\Sigma_g)$ onto 
the free group $F_n$, then $n\leq g$ by Theorem~\ref{thm:ontoFn}. It follows that $n\leq g(F_n)$.
\end{proof}

\begin{corollary} \label{cor:freeprod}
Let $m\geq 2$ be an integer and let $F_n$ denote the free group of rank $n$. Let $F_n*\Z_m$
be the free product of $F_n$ with the cyclic group of order $m$.
Then $n+1\leq g(F_{n}*\Z_m)\leq 2n+2$. 
\end{corollary}
\begin{proof} 
It can be shown that the fundamental group of 
\[ X^{2n+2}_W(1,b_1,b_2,\ldots,b_{2n+2}, a_1^m)\]
is isomorphic to $F_{n}*\Z_m$. Hence, $g(F_{n}*\Z_m)\leq 2n+2$. 

Suppose that $g(F_{n}*\Z_m)=h$, so that there is an epimorphism $\varphi:\pi_1(\Sigma_h)\to F_{n}*\Z_m$.
The kernel of the epimorphism $\phi:F_{n}*\Z_m\to \Z_m$ mapping $F_n$ to zero and a generator of 
$\Z_m$ to a generator of $\Z_m$ is a free group $F$ of rank $mn$. Since the index of $F$ in $F_{n}*\Z_m$
is $m$, so is the index of $\phi^{-1} (F)$ in $\pi_1(\Sigma_h)$. Thus $\phi^{-1} (F)$ is isomorphic 
to the fundamental group of a closed orientable surface of genus $m(h-1)+1$, an $m$-sheeted 
covering of $\Sigma_h$. By Theorem~\ref{thm:ontoFn}, $m(h-1)+1 \geq mn$. It follows now that
$h\geq n+1$.
\end{proof}

\bigskip
\subsection{Abelian groups}
Let $\Z^n$ denote the free abelian group of rank $n$ and let $\Z_m$ denote
the cyclic group of order $m$ for $m\geq 2$. 
By considering the standard presentation  
$$\langle x_1,x_2,\ldots, x_n \, \vert \, x_ix_jx_i^{-1}x_j^{-1} , 1\leq i<j\leq n \rangle$$
of the free abelian group $\Z^n$,
one can conclude from Theorem~\ref{thm:genus}~(a) that 
$g(\Z^n)\leq 3n^2-n$. By reexamining the proof of 
Theorem~\ref{thm:main} for this special case we can get a better estimate. 
In fact, we will give better bound for any finitely presented abelian group.

Let $m_i\geq 2$. For the abelian group $\Z^n\oplus \Z_{m_1} \oplus\cdots \oplus \Z_{m_k} $, 
consider the presentation
\begin{eqnarray}
\left< x_1,x_2,\ldots, x_{n+k} \; | \; x_ix_jx_i^{-1}x_j^{-1}, (x_{n+s})^{m_s} \right>,
\end{eqnarray}
where $1\leq i<j\leq n+k$ and $1\leq s\leq k$.

\begin{figure}[hbt]
 \begin{center}
    \includegraphics[width=11cm]{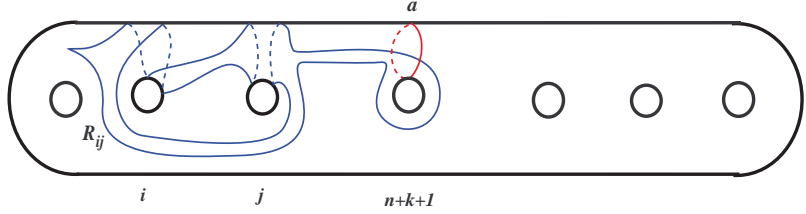}
  \caption{The curve $R_{ij}$ on the surface $\Sigma_{(2n+2k+1)}$.}
  \label{abelian}
   \end{center}
 \end{figure}

\begin{theorem} \label{thm:freeabelian}
Let $n$ and $k$ be nonnegative integers with $n+k\geq 3$ 
and let $m_i\geq 2$. Let $\Gamma=\Z^n\oplus \Z_{m_1}\oplus \cdots \oplus \Z_{m_k}$.
Then $\frac{n+k+1}{2}\leq g(\Gamma)\leq 2n+2k+1$.
\end{theorem}

\begin{proof}
The first inequality follows from Theorem~\ref{thm:genus}(b) using the fact that the minimal number
of generators for $\Gamma$ is $n+k$ and that $\Gamma$ is not a surface group.

For the second inequality, we set $g=2n+2k+1$. For $i<j$, consider the simple loop $R_{ij}$ on the
surface $\Sigma_{g}$ of genus $g$ as shown in Figure~\ref{abelian}. 
Clearly, the loop $R_{ij}$ represents the element 
\[
a_ia_j a_i^{-1}(b_i\cdots b_j)a_j^{-1} b_{n+k+1}^{-1} (b_i\cdots b_j)^{-1} 
\]
in $\pi_1(\Sigma_{g})$. Note that each $R_{ij}$ intersects
the simple closed curve $a$ transversely only once, where $a$ is the curve appearing in the word $W$
for $g=2n+2k+1$. After setting all $b_k=1$, this word reduces to $a_ia_ja_i^{-1}a_j^{-1}$.
Also, for each $s=1,\ldots,k$, the element $(a_{n+s})^{m_s}b_{n+s}^{-1}$ is represented by
a simple loop $T_s$ which intersects $B_{n+s}$ only once.  
It follows now from Theorem~\ref{thm:presentation} that the fundamental group of
\begin{eqnarray*} 
X=X^g_W(1,b_1,b_2,\ldots,b_g,T_1,\ldots,T_k,R_{ij})
\end{eqnarray*}
is isomorphic to $\Gamma$. Here, all $R_{ij}$ are included, so that
the manifold $X$ is a fiber sum of $(g+1+k+\frac{(n+k)(n+k+1)}{2})$ copies of
$X^g_W(1)$.

This shows that $g(\Gamma)\leq g$, finishing the proof of Theorem.
\end{proof}

\bigskip
\subsection{The fundamental groups of closed nonorientable surfaces}
Let $N_g$ denote be closed connected nonorientable surface of genus $g\geq 1$, that is the 
connected sum of $g$ copies of the real projective plane. 
The fundamental group $\pi_1(N_g)$ of $N_g$ has a presentation 
\[ \langle \, c_1,c_2,\ldots,c_g \, \vert \, c_1^2c_2^2\cdots c_g^2 \, \rangle. \] 

If $g=1$, then $\pi_1(N_g)$ is isomorphic
to $\Z_2$. Hence, $g(\pi_1(N_g))=2$.

Suppose that $g\geq 2$.  Since $\pi_1(N_g)$ cannot be generated by less than $g$ elements
and since it is not isomorphic to the fundamental group of a closed orientable surface, 
we get from Theorem~\ref{thm:genus}~(b) that $g(\pi_1(N_g))\geq \frac{g+1}{2}$.

On the other hand, the loop $a_1^2a_2^2\cdots a_g^2 (b_1b_2\cdots b_g)^{-1}$ can be represented by a simple closed curve $R$, which intersects the simple closed curve $B_g$ only once. It follows that
the fundamental group of $X_W^{2g} (1, b_1,\ldots,b_g, R)$ is isomorphic to  $\pi_1(N_g)$.
Therefore, $g(\pi_1(N_g))\leq 2g$.

\bigskip
\subsection{The braid group} 
\begin{figure}[hbt]
 \begin{center}
    \includegraphics[width=11cm]{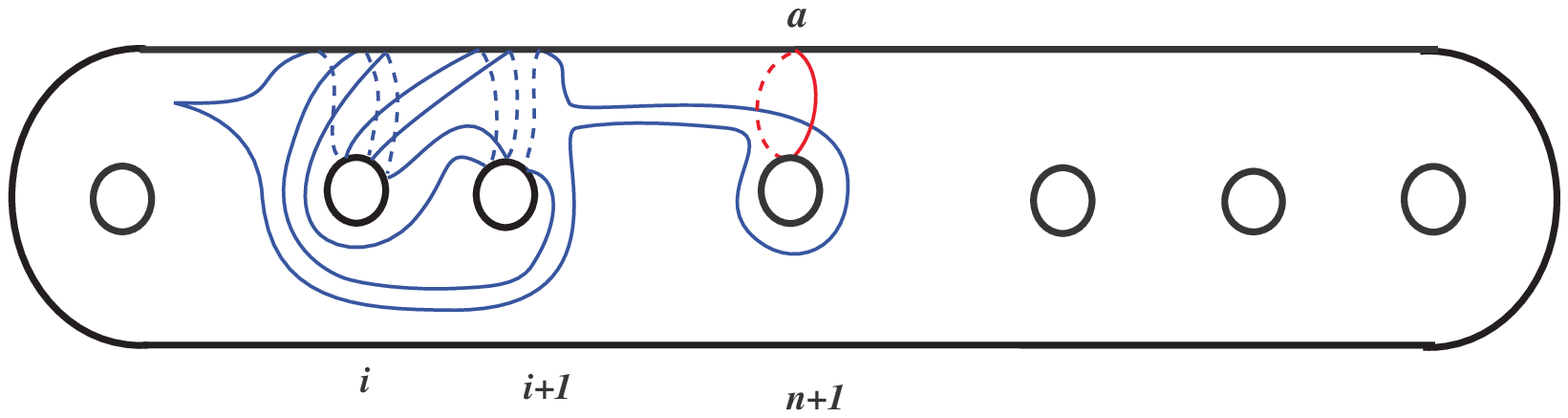}
  \caption{}
  \label{braid}
   \end{center}
 \end{figure}

The braid group $B_n$ on $n$ strands has a presentation
\begin{eqnarray*}
\left< \sigma_1,\sigma_2,\ldots, \sigma_{n-1} \; | \; \sigma_i\sigma_{i+1}\sigma_{i}\sigma_{i+1}^{-1}\sigma_i^{-1}\sigma_{i+1}^{-1},
\sigma_j\sigma_k\sigma_j^{-1}\sigma_k^{-1} \;\; |j-k|\geq 2
  \right>.
\end{eqnarray*}
On the oriented surface $\Sigma_{2n+1}$, consider the simple loops as in 
Figure~\ref{abelian} for the relations of the type $\sigma_j\sigma_k\sigma_j^{-1}\sigma_k^{-1}$
and the simple loops in Figure~\ref{braid} for the relations of the type
$\sigma_i\sigma_{i+1}\sigma_{i}\sigma_{i+1}^{-1}\sigma_i^{-1}\sigma_{i+1}^{-1}$. It can now easily
be shown that there is a Lefschetz fibration of genus $2n+1$ whose fundamental group is 
isomorphic to $B_n$. It follows that $2\leq g(B_n)\leq 2n+1$.

\bigskip

\subsection{The group $SL(2,\Z)$} The group $SL(2,\Z)$ has a presentation with two generators
$x$ and $y$ and with two relators $x^2y^3$ and $x^4$. For this presentation, 
we have $d(SL(2,\Z))=3$. Since $SL(2,\Z)$ is not a surface group, 
by Theorem~\ref{thm:genus} we get the bounds
$2\leq g(SL(2,\Z))\leq 6$. However, by repeating the proof of Theorem~\ref{thm:main}
for $SL(2,\Z)$ can see that there is a genus $4$ Lefschetz fibration whose total
space has fundamental group isomorphic to $SL(2,\Z)$. The construction of the loop $R$
corresponding to the relator $x^2y^3$ can be done without increasing the genus.
Just take $R$ to be a simple representative of the loop $a_1^2 a_2^3 b_2^{-1}b_1^{-1}$.

\bigskip

\section{Other invariants and problems} \label{section6}
In~\cite{kot}, Kotschick defines two invariants for finitely presented groups.
For a finitely presented group $\Gamma$, these invariants are defined as
\[
q(\Gamma)=\inf \{ \chi (X) \}
\]
and
\[
p(\Gamma)=\inf \{ \chi (X)-\abs{\sigma(X)} \},
\]
where infimums are taken aver all closed orientable $4$-manifolds $X$ with 
$\pi_1(X)\cong \Gamma$. Here,  $\chi (X)$ and $\sigma(X)$ denote the Euler 
characteristic and the signature of $X$.

For a finitely presented group $\Gamma$, the invariants $q(\Gamma)$ and $p(\Gamma)$ 
satisfy the following inequalities (cf.~\cite{kot}): 
$$2-2b_1(\Gamma)+b_2(\Gamma) \leq q(\Gamma) \leq 2 (1-d(\Gamma)),$$
and
$$2-2b_1(\Gamma) \leq p(\Gamma) \leq q(\Gamma).$$
It follows that
$$2-4g(\Gamma) \leq p(\Gamma) \leq q(\Gamma).$$

\bigskip

Here are some problems motivated by the discussions in this paper and the
properties of the invariants $p$ and $q$ proved in~\cite{kot}.

{\bf Problem 1.} Find the exact values of $g(F_n)$ and $g(\Z^n)$. 

{\bf Problem 2.} Compare $g(\Gamma_1\times \Gamma_2)$ and $g(\Gamma_1\! * \Gamma_2)$ 
with $g(\Gamma_1)$ and $g(\Gamma_2)$ for any two finitely
presented groups $\Gamma_1$ and $\Gamma_2$, where $\Gamma_1\! * \Gamma_2$ is 
the free product of $\Gamma_1$ with $\Gamma_2$.

{\bf Problem 3.} For a finitely presented group $\Gamma$ and a subgroup 
$\Gamma'$ of finite index $k$ in $\Gamma$, compare the values 
$g(\Gamma')$ and $g(\Gamma)$.

{\bf Problem 4.} In all examples I know, Lefschetz fibrations of genus $2$ 
having singular fibers have abelian fundamental group. 
Is this true for all genus-$2$ Lefschetz fibrations?

\end{document}